\documentclass[10pt]{amsart}
\usepackage[latin1]{inputenc}
\usepackage{amsfonts}
\usepackage{amsmath}
\usepackage{amssymb}
\usepackage{amsthm}
\usepackage{fontenc}

\usepackage{graphicx}
\usepackage[margin=1.4in]{geometry}

\newtheorem{theorem}{Theorem}[section]
\newtheorem{lemma}[theorem]{Lemma}
\newtheorem{proposition}[theorem]{Proposition}

\newtheorem{corollary}[theorem]{Corollary}

\newtheorem{conjecture}{Conjecture}

\numberwithin{equation}{section}

\DeclareMathOperator{\Avg}{Avg}
\DeclareMathOperator{\Trace}{Trace}
\DeclareMathOperator{\Mid}{Mid}
\DeclareMathOperator{\supp}{supp}
\DeclareMathOperator{\Hom}{Hom}
\DeclareMathOperator{\colim}{colim}

\title{Average-Value Tverberg Partitions via Finite Fourier Analysis} 

\author{Steven Simon}

\begin{document}

\maketitle

\begin{abstract}

The long-standing topological Tverberg conjecture claimed, for any continuous map from the boundary of an $N(q,d):=(q-1)(d+1)$-simplex to $d$-dimensional Euclidian space, the existence of $q$ pairwise disjoint subfaces whose images have non-empty $q$-fold intersection. The affine cases, true for all $q$, constitute Tverberg's famous 1966 generalization of the classical Radon's Theorem.  Although established for all prime powers in 1987 by \"Ozaydin, counterexamples to the conjecture, relying on 2014 work of Mabillard and Wagner, were first shown to exist for all non-prime-powers in 2015 by Frick. Starting with a reformulation of the topological Tverberg conjecture in terms of harmonic analysis on finite groups, we show that despite the failure of the conjecture,  continuous maps \textit{below} the tight dimension $N(q,d)$ are nonetheless guaranteed $q$ pairwise disjoint subfaces -- including when $q$ is not a prime power -- which satisfy a variety of ``average value" coincidences, the latter obtained as the vanishing of prescribed Fourier transforms. 
\end{abstract}

\section{Introduction} 

\subsection{A History of the Topological Tverberg Problem} The celebrated Tverberg Theorem [23] of 1966 states that any $(q-1)(d+1) +1$ points in $\mathbb{R}^d$ can be partitioned into $q$ pairwise disjoint sets whose convex hulls have non-empty $q$-fold intersection. The case $q=2$ is the classical Radon's Theorem of 1921 [18].  Viewing intersecting convex hulls as the images of pairwise disjoint subfaces under a piecewise-affine map from the boundary of a $(q-1)(d+1)$-dimensional simplex to $d$-dimensional Euclidian space, a long-standing conjecture and ``holy grail" [12] of topological combinatorics claimed a continuous extension:
		
 \begin{conjecture} Any continuous map $f: \partial \Delta^{(q-1)(d+1)} \rightarrow \mathbb{R}^d$ admits a topological Tverberg partition, i.e., $q$ pairwise disjoint subfaces $\sigma_1,\ldots, \sigma_q$ such that $f(\sigma_1)\cap \cdots \cap f(\sigma_q)\neq \emptyset$. \end{conjecture}
 
That $N(q,d):=(q-1)(d+1)$  is tight for given $q$ and $d$ can be easily seen by considering piecewise-affine maps in general position and counting the codimension of their intersections (see, e.g., [13, 1.1 Remark 3]). Positive solutions to the conjecture were established for $q=2$ in 1979 by Bajm\'oczy and B\'ar\'any [2], for odd primes by B\'ar\'any, Shlosman, and Sz\"ucs [3] in 1981, and for all prime powers by \"Ozaydin [17] in 1987:

\begin{theorem} [Topological Tverberg] Any continuous map $f:\partial \Delta^{(q-1)(d+1)}\rightarrow \mathbb{R}^d$ admits a topological Tverberg partition if $q$ is a prime power. \end{theorem}

An essential ingredient in \"Ozaydin's proof is the reduction of a topological Tverberg partition to the existence of zeros of an induced $S_q$-equivariant map $f^{\times q}: (\partial\Delta^{N(q,d)})^{\times q}_{(2)}\rightarrow_{S_q} W_q^d$, where for any simplicial complex $K$ the \textit{deleted $q$-fold product} \begin{equation} K^{\times q}_{(2)} = \{(x_1,\ldots, x_q)\in K^q\mid \supp(x_i)\cap \supp(x_j)=\emptyset\,\, \textit{for}\,\, i\neq j\} \end{equation} consists of all $q$-tuples of points in $K$ with pairwise disjoint support, $W_q=\{(w_1,\ldots, w_q)\in \mathbb{R}^q\mid \sum_{j=1}^q w_j=0\}$, and the symmetric group $S_q$ acts on both spaces by permuting coordinates. That any  $h: (\partial \Delta^{N(q,d)})^{\times q}_{(2)}\rightarrow_{S_{q}} W_q^d$ must vanish for $q$ a prime power [17, Corollary 3.4] yielded Theorem 1.1. That such an equivariant map exists without zeros [17 ,Theorem 4.2] for all other $q$ -- but was not shown to arise from some $f:\partial \Delta^{N(q,d)}\rightarrow \mathbb{R}^d$ -- was the main reason the conjecture remained open for so long. 

	The central breakthrough arrived in the 2014 extended abstract [14] of Mabillard and Wagner, who showed (Theorem 3) the \textit{equivalence}, for any $(q-1)k$-dimensional simplex $K$ with $k\geq 3$, between (i) equivariant maps $h: K^{\times q}_{(2)}\rightarrow_{S_q} W_q^d$ without zeros and (ii) continuous maps $f:K\rightarrow \mathbb{R}^d$ for which any $q$ pairwise disjoint faces of $K$ have images with empty $q$-fold intersection. Using this equivalence and \"Ozaydin's result [17, Lemma 4.1] on the existence of certain $S_q$-equivariant maps, Frick ([9] and again in [4]) deduced counterexamples to the generalized Van Kampen conjecture for all non-prime-powers $q$, and then, as an application of the ``constrained" method of Blagojevi\'c, Frick, and Ziegler [5], to Conjecture 1 itself even for maps from the \textit{full} simplex $\Delta^{N(q,d)}$. These initial counterexamples occurred in dimensions $d=qk+1$ with $k\geq 3$, with a subsequent lowering to $d=qk$ in the journal version [13] of [14]. 	  

\subsection{A Finite Harmonic Analysis Perspective} The objective of this paper is to show that, despite the failure of the topological Tverberg conjecture in general, continuous maps $f:\partial \Delta^n \rightarrow \mathbb{R}^d$ with $n$ \textit{below} the tight dimension $N(q,d)$ are nonetheless guaranteed collections of $q$ pairwise disjoint faces whose images satisfy a variety of other coincidence types, even when $q$ is no longer constrained to be a prime power. As discussed in Section 3, our starting point is a reformulation of the topological Tverberg problem in terms of harmonic analysis on finite abelian groups. Full Tverberg partitions are shown to be equivalent to the vanishing of all Fourier transforms except those arising from the trivial representation (Lemma 3.1), while the annihilation of fewer prescribed transforms results in ``average-value" Tverberg partitions (Theorems 2.1 and 2.2) which have natural interpretations in terms of barycenters (Section 2.1).
	
	The underlying equivariant topological methods, discussed in Section 4, are no different than those standardly employed in topological combinatorics -- here, we are ultimately reduced to calculating  characteristic classes in group cohomology for linear representations determined by the transforms considered. Nonetheless, the polynomial conditions of the central technical result Theorem 3.2 show that, by carefully selecting transforms, the Fourier perspective affords new applications of these classical techniques. This approach was first introduced by the author [21] in the context of mass partition problems, another central topic of the field, where transform annihilation produced a variety of equipartition-types, again by non-prime-power numbers of regions. We conclude the paper with some comments (Section 4.3) on the ``homotopic" optimality of Theorem 3.2 and its corollaries Theorems 2.1--2.2, a situation which can be compared to that of the topological Tverberg problem itself and the existence of $S_q$-equivariant maps.

\section{Average-Value and Barycentric Tverberg Partitions} 

Our main theorems can be summarized as follows. Let $q=rp^k$, $p$ prime, and let $0\leq a \leq k$. For appropriate $n$, we show that for any continuous map $f:\partial \Delta^n\rightarrow \mathbb{R}^d$, there exist $q$ points, one from each of $q$ disjoint faces, which can be split into $q'=rp^{k-a}$ subfamilies $\mathcal{F}_1=\{x_{1,1},\ldots x_{1, p^a}\},\ldots, \mathcal{F}_{q'}=\{x_{q',1},\ldots, x_{q',p^a}\}$ of $p^a$ points each such that \vspace{.075in} 
 
 \noindent (i) For each $1\leq j \leq q'$, the points of each $\mathcal{F}_j$ have identical image (thus $\cap_{i=1}^{p^a}f(\sigma_{j,i})\neq\emptyset$ for the corresponding collections of faces), and additionally \vspace{.075in}
 
 \noindent (ii) For the representatives $\{x_1:=x_{1,1},\ldots, x_{q'}:=x_{q',1}\}$ with $f(x_j)\in \cap_{i=1}^{p^a} f(\sigma_{j,i})$, one can guarantee in the case of Theorem 2.1 that, for varying $\ell\geq 2$, there are relatively large numbers of symmetric order $\ell$ subsets $\{x_{j_1},\ldots, x_{j_\ell}\}$ whose \textit{average values} \begin{equation} \Avg(f;x_{j_1},\ldots, x_{j_\ell}):=\frac{f(x_{j_1})+\cdots + f(x_{j_\ell})}{\ell} \end{equation} are fixed, and likewise in Theorem 2.2 that the $\{x_1,\ldots, x_{q'}\}$ can be further subdivided into $p^{k-a}$ sets, each containing $r$ points, which have the same average-value. 
    
\begin{theorem} Let $q=rp^k$, $p$ prime, and for any $0\leq a \leq k$, let $G=H\oplus \bar{G}$, where $H=\mathbb{Z}_p^{a}$ and $\bar{G}=\mathbb{Z}_p^{k-a}\oplus\mathbb{Z}_r$. Then for any continuous map $f:\partial \Delta^n \rightarrow \mathbb{R}^d$, $n=(d+1)(q-1)-d[(k-a)(p-1)+r-1]$, there exist $q$ points $\{x_g\mid x_g\in\sigma_g \}_{g\in G}$ from $q$ pairwise disjoint faces such that\vspace{.075in}

\noindent (i) $f(x_{h+\bar{g}})= f(x_{\bar{g}})$ for all $h\in H$ and all $\bar{g}\in \bar{G}$, and\\ 
(ii) For any $\ell \leq \min \{r,p\}$ and each of the $\binom{p}{\ell}^{k-a}\binom{r}{\ell}$ subsets $S=S_{a+1} \times \cdots \times S_{k+1} \subset \bar{G}$ with $|S_j|=\ell$ for all $a+1\leq j \leq k+1$, the $(\ell!)^{k-a}$ average-values \begin{equation} \Avg(f;  x_s, x_{\phi(s)}, \ldots, x_{\phi^{\ell-1}(s)})= c(S) \end{equation} are constant for all $s\in S$ and all permutations $\phi=(\phi_{a+1},\ldots, \phi_{a+k})$ for which each $\phi_j$ is an $\ell$-cycle of $S_j$. 
\end{theorem}

\begin{theorem}  Let $q=rp^k$, $p$ prime, and for any $0\leq a\leq k$, let $q'=rp^{k-a}$. Then for any continuous map $f:\partial\Delta^n \rightarrow \mathbb{R}^d$, $n=(d+1)(q-1)-d[p^{k-a}(r-1)-1]$, there exist $q$ points $x_1,\ldots, x_q$, one each from $q$ pairwise disjoint faces, such that (i) $f(x_j)= f(x_{j+i})$ for all $1\leq i \leq p^a$ and all $1\leq j \leq q'$, and (ii) \begin{equation} \Avg(f; x_1,\ldots, x_r) = \Avg (f; x_{r+1}, \ldots, x_{2r}) = \cdots = \Avg (f; x_{q'-r+1},\ldots, x_{q'}) \end{equation} 
 \end{theorem}
 
As a special case of Theorem 2.2, we note that letting $p=2$, $k=1$, and $a=0$ gives an alternating-sum generalization of the topological Radon theorem Bajm\'oczy and B\'ar\'any [2]. In the planar cases, this forces all but at most two of the disjoint faces to be vertices:

\begin{corollary} For any continuous map $f: \partial \Delta^{2r+1}\rightarrow \mathbb{R}^2$, there exist either 

\begin{itemize} \item $2r-1$ vertices $x_1,\ldots, x_{2r-1}$ and a point  $x_{2r}$ from the remaining disjoint 2-dimensional face, or 
\item $2r-2$ vertices $x_1,\ldots, x_{2r-2}$ and points $x_{2r}$ and $x_{2r+1}$, one from each of the remaining disjoint edges, for which \end{itemize} \begin{equation} \sum_{j=1}^{2r} (-1)^j f(x_j) = 0 \end{equation} 
\end{corollary}

\subsection{A Barycentric Interpretation} The average-value conditions in both theorems are vacuous when $a=k$. Theorem 1.1 is then recovered when $r=1$, while for $r>1$ one has its $r$-fold iteration. From the viewpoint of Conjecture 1.1, the ``worst case" $0\leq k<a$ scenarios occur when (using the notation of the paragraph preceding Theorems 2.1--2.2) the $q'$ images $f(x_j)\in \cap_{i=1}^{p^a} f(\sigma_{j,i})$ are all distinct, lest even more faces have overlapping image than guaranteed by (i). This  situation is ``homotopically" generic, as discussed in Section 4.3. The average-value conditions (ii) then have natural geometric interpretations: in Theorem 2.1, that for each order $\ell^{k-a+1}$ subset $S\subset\bar{G}$ as above, the $(\ell !)^{k-a}$ subsets $\{f(x_s),\ldots, f(x_{\phi^{\ell-1}(s)})\}\subset \{f(x_s) \in \cap_{h\in H} f(\sigma_{h+s})\}_{s\in S}$ all have the same barycenter, and in Theorem 2.2 that the $f(x_j)\in\cap_{i=0}^{p^a}f(\sigma_{j+i})$ can be split into $p^{k-a}$ sets of $r$ points each, the barycenters of which are all identical.

	Relatively simple examples for non-prime-powers can be seen when $q=2p^2$, $p$ an odd prime, and $a=1$. If $f:\partial\Delta^n \rightarrow \mathbb{R}^d$ and $n$ is either exactly (a) $dp$ or (b) $d(p-1)$ below $N(2p^2,d)=(d+1)(2p^2-1)$, then at a minimum the $2p$ images $f(x_j)\in \cap_{i=0}^{p-1}f(\sigma_{j+2pi})$ in case (a) determine $\binom{p}{2}$ pairs of edges with equal centers: $\Mid[f(x_{j_1}), f(x_{j_2+p})]=\Mid[f(x_{j_2}),f(x_{j_1+p})]$ for any subset $\{j_1,j_2\}\subset \{0,\ldots, p-1\}$, and in case (b) can be partitioned into $p$ edges $[f(x_0),f(x_1)], \ldots, [f(x_{2p-2}), f(x_{2p-1})]$ with common midpoint.

\section{A Finite Fourier Approach} 

We now describe a finite harmonic analysis reformulation of the topological Tverberg problem, of which Theorems 2.1 and 2.2 are special cases. For a given simplex $\Delta^n$, one can index all the collections $Q=\{x_g\}_{g\in G}$ of $q$ points of $\partial \Delta^n$ with pairwise disjoint support by any fixed group of order $G$. Given any continuous map $f=(f_1,\ldots, f_d): \partial \Delta^n\rightarrow \mathbb{C}^d$, each collection thereby determines $d$ maps $F_1,\ldots, F_d: G\rightarrow \mathbb{C}$, $g\mapsto f_i(x_g)$, each of which has Fourier expansion $F_i(g)=\sum_{\chi \in \hat{G}} d_\chi \Trace(c_{i,\chi}\,\chi_g)$, where $\hat{G}$ consists of the distinct irreducible unitary representations $\chi: G\rightarrow U(d_\chi)$ and the $c_{i,\chi}=\frac{1}{|G|}\sum_{u\in G} F_i(u)\chi_u^{-1}\in M(\mathbb{C},d_\chi)$ are the corresponding matrix-valued Fourier transforms (see, e.g., [22, Theorem 5.5.4]). Choosing $G=\mathbb{Z}_{q_1}\oplus \cdots \oplus \mathbb{Z}_{q_k}$ to be abelian, the irreducible representations are all one-dimensional and naturally indexed by the group itself (see, e.g., [20, Theorem 9]), so one has the particularly simple decomposition \begin{equation} F_i(g) = \sum_{\epsilon \in \oplus_{j=1}^k\mathbb{Z}_{q_j}} c_{i,\epsilon} \chi_\epsilon (g), \end{equation} 
where
\noindent \begin{equation} c_{i,\epsilon} = \frac{1}{|G|} \sum_{u \in G} f_i(x_u) \chi_\epsilon^{-1}(u) \in \mathbb{C} \end{equation} and the $\chi_\epsilon: G\rightarrow U(1)$ are given explicitly by $\chi_\epsilon(g)=\Pi_{j=1}^k \zeta_{q_j}^{\epsilon_jb_j}$ for each $\epsilon=(\epsilon_1,\ldots, \epsilon_k)$ and each $g=(b_1,\ldots, b_k)\in G$, $\zeta_{q_j}=\exp(2\pi i/q_j)$.\vspace*{.075in}

As a preliminary observation, we note that topological Tverberg partitions are equivalent to the annihilation of all $(q-1)d$ transforms not arising from the trivial representation:

\begin{lemma} A continuous map $f:\partial \Delta^n\rightarrow \mathbb{R}^d$ admits a full Tverberg partition iff there exists some $\{x_g\in \sigma_g\}$ from disjoint $\sigma_g$ for which $c_{i,\epsilon}=0$ for all $1\leq i\leq d$ and all $\epsilon\neq 0$. 
\end{lemma}

\begin{proof} As $G$ is abelian, the $\chi_\epsilon$ form a basis for the complex vector space consisting of all maps $h: G\rightarrow \mathbb{C}$ (and in fact an orthonormal basis with respect to the inner product $\langle h_1, h_2\rangle = \frac{1}{|G|} \sum_{u\in G} h_1(u)\overline{h_2(u)}$, see, e.g. [20, Theorem 6]). On the other hand, a full Tverberg partition for $f:\partial \Delta^n\rightarrow \mathbb{R}^d$ is equivalent to each $F_i =c_{i,0} + \sum_{\epsilon\neq 0}c_{i,\epsilon} \chi_\epsilon$ being the constant map.\end{proof} 
	
	The following theorem gives general conditions for which the vanishing of transforms can be ensured.  Note that $\overline{c_{i,\epsilon}}=c_{i,-\epsilon}$ when each $f_i$ above is real-valued, and that $\chi_\epsilon$ is real-valued iff $\epsilon$ has order 2.
	
\begin{theorem} Let $q=q_1\cdots q_k$, and let $\epsilon_1, \ldots, \epsilon_m \in \mathbb{Z}_{q_1}\oplus\cdots\oplus\mathbb{Z}_{q_k}$, $\epsilon_j=(\epsilon_{j,1},\ldots, \epsilon_{j,k})$.\\
(a) Let $n=2dm+q-1$. If \begin{equation} h(y_1,\ldots, y_k)= \Pi_{j=1}^m (\epsilon_{j,1}y_1+\cdots +\epsilon_{j,k}y_k)^d \end{equation}  is non-zero in $\mathbb{Z}[y_1,\ldots, y_k]/(q_1y_1,\ldots, q_ky_k)$, then for any continuous map $f: \partial \Delta^n \rightarrow \mathbb{C}^d$, there exist $q$ points of $\partial \Delta^n$ with pairwise disjoint support such that $c_{i,\epsilon_j}=0$ for each  $1\leq j \leq m$ and each $1\leq i \leq d$ in the Fourier expansion (3.1).\\
\noindent (b) Suppose that $d$ is odd, $q_1=2r_1,\ldots, q_{k'}=2r_{k'}$ are even, that $\epsilon_1 ,\ldots, \epsilon_{m'}$ are the elements of order 2, and let $n=d(2m-m')+q-1$. If 
\begin{equation} h(x_1,y_1,\ldots, x_{k'}, y_{k'}) =\Pi_{j=1}^{m'} (\epsilon_{j,1}x_1+\cdots +\epsilon_{j,k'}x_{k'})^d \Pi_{j=m'+1}^m (\epsilon_{j,1}y_1 +\cdots +\epsilon_{j,k'}y_{k'})^d\end{equation} is non-zero in $\mathbb{Z}_2[x_1,y_1,\ldots, x_{k'},y_{k'}]/(x_1^2-r_1y_1,\ldots, x_{k'}^2-r_{k'}y_{k'})$, then for any continuous map $f: \partial \Delta^n \rightarrow \mathbb{R}^d$, there exist $q$ points of $\partial \Delta^n$ with pairwise disjoint support such that $c_{i,\epsilon_j}=c_{i,-\epsilon_j}=0$ for each  $1\leq j \leq m$ and each $1\leq i \leq d$ in the Fourier expansion (3.1).
 \end{theorem} 
 
\subsection{Proof of Theorems 2.1 and 2.2} We defer the proof of Theorem 3.2 to Section 4. Note that for $G=\mathbb{Z}_p^k\oplus \mathbb{Z}_r$, $p$ prime, the polynomials (3.3) and (3.4) are non-zero provided each $\epsilon_j \notin 0\oplus\mathbb{Z}_r$: (3.3) is non-zero in $\mathbb{Z}_p[y_1,\ldots, y_k]$ even after quotienting by $(y_{k+1})$, and likewise (3.4) is non-zero in $\mathbb{Z}_2[x_1,\ldots, x_k]$ after quotienting by $(x_{k+1},y_{k+1})$. Theorems 2.1 and 2.2 follow easily:

\begin{proof} For Theorem 2.1, let $f:\partial \Delta^n\rightarrow \mathbb{R}^d$, $n=(d+1)(q-1)-d[(k-a)(p-1)+r-1]$. We have $G=H\oplus \bar{G}$, where $H=\mathbb{Z}_p^a$ and $\bar{G}=\mathbb{Z}_p^{k-a}\oplus\mathbb{Z}_r$. For each $1\leq i \leq d$, we annihilate the $m=q-[(k-a)(p-1)+r]$ coefficients $c_{i,\epsilon}$ with $\epsilon \notin \cup_{j=a+1}^k \mathbb{Z}_{p}\mathbf{e}_j\cup \mathbb{Z}_r\mathbf{e}_{k+1}$, where $\mathbf{e}_j$ is the $j$-th standard basis vector of $G$. The polynomial conditions are met by the observation above, and that the dimension conditions are satisfied is verified case-by-case: $\mathbb{R}^d=\mathbb{C}^{d/2}$ if $d$ is even, and the annihilation of the resulting transforms of $f:\partial \Delta^n \rightarrow \mathbb{C}^{d/2}$ gives $n= 2(d/2)m+q-1$, so part (a) of Theorem 3.2 applies. If both $d$ and $p$ are odd, then none of the $\epsilon$ have order 2, and as $c_{i,-\epsilon}=\overline{c_{i,\epsilon}}$ and $m$ must be even, one only need annihilate half of these $c_{i,\epsilon}$ (e.g., those $\epsilon$ with first non-zero coordinate $\epsilon_{j_0}\leq \frac{p-1}{2}$). Hence $n=2d(m/2)+q-1$, and again one may apply part (a). Finally, if $d$ is odd and $p=2$, then $m' =2^k-(k-a+1)$ if $r$ is odd and $m'=2^{k+1}-(k-a+2)$ if $r$ is even. Annihilating half of the remaining (necessarily even) coefficients gives $n=d(2m-m')+q-1$, so part (b) may be used. Writing each $g\in G$ as $g=h+\bar{g}$, $h=(b_1,\ldots, b_a)\in H$ and $\bar{g}=(b_{a+1},\ldots b_{k+1})\in \bar{G}$, in all cases the resulting Fourier expansion is
	\begin{equation} F_i(g)= c_{i,0} + \sum_{u=1}^{p-1}c_{i, u\mathbf{e}_{a+1}}\zeta_p^{ub_{a+1}} +\cdots +\sum_{u=1}^{r-1}c_{i,u\mathbf{e}_{k+1}}\zeta_r^{ub_{k+1}} \end{equation}
Thus (i) $F_i(h+\bar{g})=F_i(\bar{g})$ for all $h\in H$ and all $\bar{g}\in \bar{G}$, and (ii) if $\ell\leq \min\{p,r\}$ and $S=S_{a+1} \times \cdots \times S_{a+k} \subset \bar{G}$ with $|S_j|=\ell$ for all $a+1\leq j \leq k+1$, it follows that the sums $F_i(s) + F_i(\phi(s)) + \cdots + F_i(\phi^{\ell-1}(s))$ are constant for all $s\in S$ and all $k$-tuples $\phi=(\phi_{a+1},\ldots, \phi_{k+1})$ of $\ell$-cycles of the $S_j$. 

	For Theorem 2.2, we again let $G=H\oplus \bar{G}$, but now we annihilate all $c_{i,\epsilon}$ with $\epsilon \in \mathbb{Z}_p^{k-a}\oplus 0\subset \bar{G}$ in addition to all $\epsilon\notin 0\oplus \bar{G}$. The same arguments as above show that the dimensional and polynomial conditions of Theorem 3.2 are satisfied, and as before $F_i(h+\bar{g})=F_i(\bar{g})$ for all $h\in H$ and all $\bar{g}\in \bar{G}$. Letting $H'=\mathbb{Z}_p^{k-a}$, one now has $F_i(\bar{g})=F_i(h',b)=c_{i,0} + \sum_{\epsilon\in \bar{G}, \epsilon_{k+1}\neq 0}c_{i,\epsilon}\chi_\epsilon(h',0)\zeta_r^{\epsilon_{k+1}b}$ for all $h'\in H'$ and all $b\in \mathbb{Z}_r$, and therefore that $F_i(h',0)+\cdots +F_i(h',r-1)=rc_{i,0}$ is constant. \end{proof}

\section{Topological Underpinnings}

	Our proof of Theorem 3.2 follows the usual \textit{configuration-space/test-map scheme}, the standard method for the reduction of problems in discrete and combinatorial geometry to corresponding ones of algebraic topology. See, e.g., [15, 25--26] for introductions. 
		
\subsection{Configuration-Spaces and Equivariant Test Maps} As in [3, 17, 24], all collections of $q$ points of $\partial \Delta^n$ with pairwise disjoint support can be parametrized by the deleted $q$-fold product \begin{equation} X:= (\partial \Delta^n)^{\times q}_{(2)}=\{x=(x_1,\ldots, x_q) \in \sigma_1\times\cdots\times\sigma_q\mid \sigma_i\cap\sigma_j=\emptyset \,\, \textit{for}\,\, i\neq j\}, \end{equation} which is therefore the natural \textit{configuration space} for the problem. Although the full symmetric group $S_q$ acts freely on this simplicial complex by permuting coordinates,  we shall (as in [19], and similarly in [24]) restrict to the free $G$-action induced from left multiplication after indexing $\{1,\ldots, q\}$ by the given abelian group $G=\mathbb{Z}_{q_1}\oplus \cdots \mathbb{Z}_{q_k}$. 

	For $n$ as in Theorem 3.2 and a given map $f$, evaluating the various Fourier transforms $c_{i,\epsilon}$ produces a continuous \textit{test map} $\mathcal{F}: X\rightarrow \mathbb{C}^{dm}$ in part (a) and $\mathcal{F}: X\rightarrow \mathbb{R}^{dm'}\oplus \mathbb{C}^{d(m-m')}=\mathbb{R}^{d(2m-m')}$ in part (b): 
			\begin{equation} \mathcal{F}: x \mapsto \frac{1}{|G|} \sum_{g \in G} f_i(x_g) \chi_{\epsilon_j}^{-1}(g), \end{equation}
$1\leq i \leq d$ and $1 \leq j \leq m$. The group $G$ acts linearly on the respective target spaces $\mathbb{C}^{dm}$ and $\mathbb{R}^{d(2m-m')}$ via the representation \begin{equation} \rho=\oplus_{j=1}^m\chi_{\epsilon_j}^{\oplus d} \end{equation} determined by the given transforms.\\

For convenience, throughout the remainder of this section we shall let $\mathbb{K}=\mathbb{C}$, $\ell=dm$, and $n=2\ell+q-1$ for part (a) of Theorem 3.2,  and likewise $\mathbb{K}=\mathbb{R}$, $\ell=d(2m-m')$, and $n=\ell+q-1$ for part (b). The existence of the desired collection $\{x_g\}_{g\in G}$ of points with pairwise disjoint support and prescribed vanishing coefficients in (3.1) is equivalent to a zero of the map $\mathcal{F}: X\rightarrow \mathbb{K}^\ell$. Crucially, the formulas for the Fourier transforms show that this map is automatically equivariant with respect to the two actions considered, so that a zero can be guaranteed if it can be shown more generally that \textit{any} such equivariant map $h: X\rightarrow_G \mathbb{K}^{\ell}$ vanishes given the assumptions on the polynomials (3.3) and (3.4), respectively. This is the main content of the proposition below.

\begin{proposition} Let $G=\mathbb{Z}_{q_1}\oplus\cdots \oplus\mathbb{Z}_{q_j}$ act on $X:=(\partial \Delta^n)^{\times q}_{(2)}$ via the action described above and on $\mathbb{K}^{\ell}$ via the representation (4.3),  where $\ell=dm$ and $n=2\ell+q-1$ if $\mathbb{K}=\mathbb{C}$, and $\ell=d(2m-m')$ and $n=\ell+q-1$ if $\mathbb{K}=\mathbb{R}$. \vspace{.075in}

\noindent (a) Any equivariant map $h:X \rightarrow_G \mathbb{C}^\ell$ vanishes iff the polynomial (3.3) is non-zero.\\ 
(b) For $d$ odd, any equivariant map $h: X\rightarrow_G\mathbb{R}^\ell$ vanishes if the polynomial (3.4) is non-zero.\\
(c) For $d$ odd, the vector bundle $E=X\times_G\mathbb{R}^\ell$ (4.4) below is orientable iff $\epsilon_{1,i}+\cdots +\epsilon_{m',i}=0$ for each $1\leq i \leq k'$. 
\end{proposition}

\begin{proof} The proof follows the usual constructions using the theory of vector bundles and characteristics classes (see, e.g., the standard references [11, 16]). Quotienting $X\times \mathbb{K}^{\ell}$ by the diagonal $G$-action produces the vector bundle \begin{equation} \mathbb{K}^{\ell} \hookrightarrow E:=X \times_G \mathbb{K}^{\ell} \rightarrow \overline{X}:=X/G \end{equation} For a given equivariant map $h: X \rightarrow_G\mathbb{K}^\ell$, the section $x\mapsto (x, h(x))$ of the trivial bundle $X\times\mathbb{K}^\ell$ thereby induces a section $s: \overline{X}\rightarrow E$ of (4.4). It is a basic fact (see, e.g., [8, Propositions I.7.2 and I.7.3]) that equivariant maps $h: X \rightarrow_G \mathbb{K}^\ell$ without zeros are equivalent to no-where vanishing sections of $E$, and moreover that such sections are precluded by showing that the top characteristic class of the bundle is non-zero. In part (a), this is the Chern class $c_{\ell}(E)\in H^{2\ell}(\overline{X};\mathbb{Z})$, and in part (b) the Stiefel-Whitney class $w_{\ell}(E)\in H^{\ell}(\overline{X};\mathbb{Z}_2)$. Thus the vanishing of equivariant maps of the proposition follows immediately from the identification given in Sections 4.2.1 and 4.2.2 of these classes with the respective polynomials (3.3) and (3.4). Additionally, the identification of $c_\ell(E)$ with (3.3) also demonstrates the ``only if" of part (a): $\dim(\overline{X})=n-q+1=2\ell$ (see below), so the vanishing of $c_\ell(E)$ (as the Euler class of the underlying $2\ell$-dimensional oriented real bundle) means that the primary and only obstruction class to a non-vanishing section is zero  (see, e.g., [7, Chapter 7.10--7.11]). Finally, part (c) follows by recalling that a real vector bundle is orientable iff its first Stiefel-Whitney class is zero, which we show below is true iff $\epsilon_{1,i}+\cdots +\epsilon_{m',i}=0$ for all $1\leq i \leq k'$. \end{proof}

\subsection{Characteristic Class Computations} 

Our desired identifications of $c_\ell(E)$, $w_\ell(E)$, and $w_1(E)$ are established by computations in group cohomology $H^*(BG; R)$ via a classical ``factorization trick." Here $R$ is any coefficient ring, and $BG$ is the homotopically unique classifying space of the group $G$, i.e, the base space of the universal bundle $G\hookrightarrow EG\rightarrow BG=EG/G$, where $EG$ is a contractible space on which $G$ acts freely (see, e.g., [11, Chapters 4.10--4.13] for a standard reference). That passing to group cohomology is permissible follows from the key technical fact, first established in [3, Lemma 1], that $(\partial \Delta^n)^{\times q}_{(2)}$ is a $(n-q)$-connected, $(n-q+1)$-dimensional CW complex. As the $G$-action on each $(\partial \Delta^n)^{\times q}_{(2)}$ is free, $EG=\colim_{n\to\infty}  (\partial \Delta^n)^{\times q}_{(2)}$ is a model for the total space of the universal bundle, with each $(\partial \Delta^n)^{\times q}_{(2)}=E_{n-q+1}G$ the $(n-q+1)$-skeleton of $EG$. Each quotient $\overline{X}=(\partial \Delta^n)^{\times q}_{(2)}/G=B_{n-q+1}G$ is therefore the $(n-q+1)$-skeleton of the classifying space $BG=\colim_{n\to\infty} B_{n-q+1}G$. 	

	With this viewpoint, $E$ is the pullback under the inclusion $i: \overline{X}\hookrightarrow BG$ of the bundle  $E_\rho$, where for any representation $\chi: G\rightarrow \mathbb{K}^t$, $E_\chi$ denotes the bundle   \begin{equation} \mathbb{K}^t \hookrightarrow  E_\chi:=EG\times_G \mathbb{K}^t\rightarrow BG, \end{equation}  whose Chern and Stiefel-Whitney classes are commonly denoted by $c_u(\chi)$ and $w_u(\chi)$ (see, e.g., [1, Appendix]). By naturality, the total Chern class of $E$ is therefore $c(E)=i^*(c(\rho))$, and likewise the total Stiefel-Whitney class is $w(E)=i^*(w(\rho))$. Cellular cohomology shows that the induced map $i^*: H^*(BG;R) \rightarrow H^*(\overline{X};R)$ is injective in all dimensions $d'\leq n-q+1$ for any choice of coefficient ring, so it suffices to show for parts (a) and (b) that $c_\ell(\rho)\in H^*(BG;\mathbb{Z})$ and $w_\ell(\rho)\in H^*(BG;\mathbb{Z}_2)$ are precisely the polynomials claimed, respectively, and likewise for part (c) that $w_1(\rho)=0$ when $d$ is odd. These are essentially classical exercises, though we provide sketches in each case for the sake of completeness.\\
	
	To begin, recall (see, e.g., [1, Appendix]) that for any finite group $G$, evaluation of the first Chern class $c_1(\tau)$ of a given representation $\tau: G\rightarrow U(1)$ gives an isomorphism \begin{equation} c_1: \Hom(G,U(1))\cong H^2(BG;\mathbb{Z}), \end{equation} where $Hom(G,U(1))$ is a group under tensor product. Thus the isomorphism can be written as $c_1(\tau_1\otimes \tau_2)=c_1(\tau_1)+c_1(\tau_2)$. One has the analogous isomorphism for real representations and Stiefel-Whitney classes: \begin{equation} w_1: \Hom(G,O(1))\cong H^1(BG;\mathbb{Z}_2)\end{equation} 
	
\subsubsection{$H^*(BG;\mathbb{Z})$.} It is is a basic fact that $H^*(B\mathbb{Z}_q;\mathbb{Z})=\mathbb{Z}[y]/(qy)$, $|y|=2$ (e.g., by identifying $B\mathbb{Z}_q$ with the infinite dimensional Lens Space $L^\infty(q)$ as in [10, Example 3.41]), and it can be seen in a number of ways (e.g., by (4.6)), that $y=c_1(\chi_1)$ may be taken to be the first Chern class of the standard representation $\chi_1:\mathbb{Z}_q\hookrightarrow U(1)$. Together with the general K\"unneth formula applied to $BG=B\mathbb{Z}_{q_1}\times\cdots \times B\mathbb{Z}_{q_k}$ and the isomorphism (4.6), this implies that $c_1(\chi_{\epsilon_j})=\epsilon_{j,1}y_1+\cdots + \epsilon_{j,k}y_k\in \mathbb{Z}[y_1,\ldots, y_k]/(q_1y_1,\ldots, q_ky_k)$ for each $\chi_{\epsilon_j}: G\rightarrow \mathbb{C}$. By the Whitney sum formula, $c(\rho)=c(\oplus_{j=1}^m\chi_{\epsilon_j}^{\oplus d})=\Pi_{j=1}^m[c(\chi_{\epsilon_j})]^d$, so one has the desired identification of $c_\ell(\rho)=\Pi_{j=1}^m (\epsilon_{j,1}y_1+\cdots + \epsilon_{j,k}y_k)^d$ with the polynomial $h(y_1,\ldots, y_k)$ from (3.3). 
		
\subsubsection{$H^*(BG;\mathbb{Z}_2)$.} It is easily seen from the cellular (co-)chain complex for $L^\infty(q)$ as in [10, Example 2.43] that $H^*(B\mathbb{Z}_q;\mathbb{Z}_2)=\mathbb{Z}_2$ when $q$ is odd. For $q=2r$ even, we first recall that $H^*(B\mathbb{Z}_2;\mathbb{Z}_2)=H^*(\mathbb{R}P^\infty;\mathbb{Z}_2)=\mathbb{Z}_2[x]$, $|x|=1$, while $H^*(B\mathbb{Z}_{2^a};\mathbb{Z}_2)=\mathbb{Z}_2[x,y]/(x^2)$ for $a>1$, where $|x|=1$ and $|y|=2$ (see, e.g., [6, Proposition 4.5.1]). Thus $H^*(B\mathbb{Z}_{2r};\mathbb{Z}_2)=\mathbb{Z}_2[x,y]/(x^2-ry)$ by the K\"unneth Formula with field coefficients. As above, the isomorphism (4.7) shows that $x=w_1(\chi_r)$, where $\chi_r:\mathbb{Z}_{2r}\rightarrow O(1)$ is the unique real 1-dimensional representation, while that $y=w_2(\chi_1)$ follows from the the discussion of $H^*(B\mathbb{Z}_q;\mathbb{Z})$ and that $w_2(\chi_1)$ is the first Chern class $c_1(\chi_1)$ reduced mod 2. Thus $H^*(BG;\mathbb{Z}_2)=\mathbb{Z}_2[x_1,y_1,\ldots, x_{k'},y_{k'}]/(x_1^2-r_1y_1,\ldots, x_{k'}^2-r_{k'}y_{k'})$, again by the K\"unneth formula.  
	
	Now we compute the Stiefel-Whitney classes. For $1\leq j \leq m'$, it follows from (4.7) that $w_1(\chi_{\epsilon_j})=\epsilon_{j,1}x_1+\cdots+\epsilon_{j,k'}x_{k'}$. On the other hand, each $E_{\chi_{\epsilon_j}}$ is complex when $m'<j \leq m$, so $w_1(\chi_{\epsilon_j})=0$ and $w_2(\chi_{\epsilon_j})=\epsilon_{j,1}y_1 +\cdots+\epsilon_{j,k'}y_{k'}$ is the mod 2 reduction of $c_1(\chi_{\epsilon_j})$. Thus $w(\rho)=\Pi_{j=1}^{m} [w(\chi_{\epsilon_j})]^d$ again by the Whitney sum formula, so $w_\ell(\rho)= \Pi_{j=1}^{m'} (\epsilon_{j,1}x_1+\cdots +\epsilon_{j,k'}x_{k'})^d \Pi_{j=m'+1}^m (\epsilon_{j,1}y_1 +\cdots +\epsilon_{j,k'}y_{k'})^d$ is the polynomial $h(x_1,y_1,\ldots, x_{k'}, y_{k'})$ from (3.4). Finally, $w_1(\rho)=d\sum_{j=1}^{m'} w_1(\chi_{\epsilon_j})$, which for $d$ odd is zero iff $\epsilon_{1,i}+\cdots +\epsilon_{m',i}=0$ for each $1\leq i \leq k'$. 
	
\subsection{Concluding Remarks} We close the paper with a discussion on the optimality of Theorem 3.2. First, note that the dimensions $n=2dm+q-1$ and $n=d(2m-m') +q-1$ in Proposition 4.1 are minimal. Indeed (see, e.g., [8, Proposition II.3.15]), when $n'<n$ there is no obstruction to extending \textit{any} non-vanishing equivariant map defined on the $0$-skeleton of $(\partial \Delta^{n'})_{(2)}^{\times q}$ to an equivariant map on all of $(\partial \Delta^{n'})_{(2)}^{\times q}$, simply on dimensional grounds. Thus one cannot hope to eliminate more than $m$ transforms in Theorem 3.2. In particular, this shows that the barycentric interpretations of Theorems 2.1--2.2 of Section 2.1 can be seen as generic in a homotopic sense, since the intersection of more of the $f(\sigma_{i,j})$ than guaranteed by condition (i) of those theorems is equivalent to the annihilation of too many transforms in the given dimension. 

	Part (a) of Proposition 4.1 likewise shows that the \textit{choice} of transforms in Theorem 3.2 is also optimal, again in that non-vanishing equivariant maps in the tight dimension will exist if the polynomial (3.3) determined by the selected transforms vanishes. This should be compared to the analogous situation for the topological Tverberg conjecture itself and the existence of $S_q$-equivariant maps without zeros iff $q$ is not a prime power [17, Theorem 4.2]. Finally, as the vector bundle $E$ of part (c) of Proposition 4.1 is often non-orientable in odd dimensions $d$, we remark that the mod 2 computations in Theorem 3.2(b) are necessary.

\section{Acknowledgements}

This research was funded by ERC advanced grant 32094 during the author's visits to the Hebrew University of Jerusalem under the supervision of Gil Kalai, whom the author thanks for support and many helpful discussions. The author is likewise grateful to the suggestions of the anonymous reviewers, which greatly improved the exposition of this manuscript, and also to Pavle Blagojevi\'c and Alfredo Hubard for useful comments. 

\bibliographystyle{plain}

\end{document}